\documentclass[11pt]{amsart}

\begin{document}
\include{cdmstdcmds}
\title[A generating function for $k$-alternating permutations]{An analytic derivation of a generating function for $k$-alternating permutations}

\author{Jean-Christophe Pain}
\address{CEA, DAM, DIF, F-91297 Arpajon, France}
\address{Universit\'e Paris-Saclay, CEA, Laboratoire Mati\`ere en Conditions Extr\^emes, F-91680 Bruy\`eres-le-Ch\^atel, France }
\email{jean-christophe.pain@cea.fr}

%


\thanks{}

\begin{abstract}
We study the inversion enumerator of permutations whose descent set is fixed to be the set of multiples of a fixed integer $k \ge 2$. For each $n$, let $\mathfrak{S}_n^{(k)}$ denote the set of permutations of $\{1,\dots,n\}$ whose descent set is exactly $\{i : k \mid i\}$, and define the polynomial
\[
a_n^{(k)}(q)=\sum_{\sigma \in \mathfrak{S}_n^{(k)}} q^{\mathrm{inv}(\sigma)}.
\]
We prove that the associated $q$-exponential generating function
\[
F_k(t;q)=\sum_{n\ge 0}\frac{t^n}{[n]_q!}\,a_n^{(k)}(q),
\]
where $[n]_q!$ denotes the $q$-factorial, admits an explicit closed form as a ratio of two $k$-periodically truncated $q$-exponential series. The proof is purely analytic and is based on a functional equation satisfied by $F_k(t;q)$, obtained via a decomposition of the $q$-exponential series into residue classes modulo $k$. Coefficient extraction yields a convolution identity involving Gaussian binomial coefficients, which uniquely determines the inversion enumerator. This provides an analytic alternative to classical inclusion--exclusion and structural combinatorial arguments for permutation classes with periodic descent constraints.
\end{abstract}

\maketitle

\section{Introduction}

Alternating permutations form a classical object in enumerative combinatorics \cite{Comtet1974,Stanley1989,Stanley2012}. A permutation $\sigma \in \mathfrak{S}_n$ is called alternating if it satisfies
\[
\sigma_1 < \sigma_2 > \sigma_3 < \sigma_4 > \cdots,
\]
and their enumeration goes back to the pioneering work of Andr\'e \cite{Andre1879}, who showed that their exponential generating function is given in terms of the tangent and secant functions. These numbers, now known as Euler numbers, have since been extensively studied; see, for instance, Ref. \cite{Foata1970} for a classical account. More generally, permutations with periodic descent patterns, such, naturally extend this framework and lead to rich enumerative structures \cite{Gessel1978,Gessel1984}. 
It is well known (see, e.g., Stanley~\cite{Stanley1989,Stanley2012}) that inversion enumerators of permutations admit natural $q$-analogues and satisfy convolution-type recurrences under structural restrictions. Define a permutation $\sigma$ in $\mathfrak{S}_n$ to be $k$-alternating when $\sigma_i>\sigma_{i+1}$ if and only if $k$ divides $i$. Using this terminology the alternating permutations are $2-$alternating.

This work originates from a problem posed by Mendes~\cite{Mendes2007} concerning the inversion enumerator of $k$-alternating permutations and its associated generating function. For \(n\ge1\), define the \(q\)-integer
\[
[n]_q=\frac{1-q^n}{1-q},
\]
and the \(q\)-factorial
\[
[n]_q! = \prod_{i=1}^n [i]_q.
\]
We want to prove that
\begin{equation}\label{resu}
\sum_{n=0}^{\infty}\frac{t^{n}}{[n]!_q}\sum_{\sigma \in \mathfrak{S}_{n}^{(k)}}q^{\operatorname{inv}(\sigma)}=\frac{1+\sum_{j=1}^{k-1}\sum_{n=0}^{\infty}(-1)^{n}\displaystyle\frac{t^{kn+j}}{[kn+j]!_q}}{\sum_{n=0}^{\infty}(-1)^{n}\displaystyle\frac{t^{kn}}{[kn]!_q}},
\end{equation}
where $\mathfrak{S}_n^{(k)}$ denotes the set of permutations of $\{1,\dots,n\}$ whose descent set is exactly
\[
\mathrm{Des}(\sigma)=\{ i\in\{1,\dots,n-1\} : k \mid i \}.
\]

The paper is organized as follows. In Section~2, we introduce the analytic framework underlying our approach and derive a functional equation satisfied by the $q$-exponential generating function associated with $k$-alternating permutations. This equation is obtained by decomposing the $q$-exponential series into residue classes modulo $k$ and leads to a simple algebraic relation between two auxiliary series. We then extract coefficients from this identity to obtain a convolution-type recurrence involving $q$-binomial coefficients, and prove that this recurrence uniquely determines the inversion enumerator. Finally, we verify that the latter coincides with the generating function of $k$-alternating permutations, thereby establishing the announced formula. The paper concludes with a brief discussion of the method and possible extensions.

\section{Analytic proof via a functional equation}

Let us define
\[
a_n^{(k)}(q) = \sum_{\sigma \in \mathfrak{S}_n^{(k)}} q^{\mathrm{inv}(\sigma)},
\]
and the associated $q$-exponential generating functio
\[
F_k(t;q) = \sum_{n \ge 0} \frac{t^n}{[n]!_q} a_n^{(k)}(q).
\]
We shall show that the series $F_k(t;q)$ defined by the $q$-exponential generating function coincides with the unique formal solution of a functional equation derived below. Let
\[
E_q(t) = \sum_{n \ge 0} \frac{t^n}{[n]!_q}
\]
be the $q$-exponential function, and decompose it into residue classes modulo $k$:
\[
E_q(t) = \sum_{j=0}^{k-1} E_j(t),
\]
where
\[
E_j(t) = \sum_{n \ge 0} \frac{t^{kn+j}}{[kn+j]!_q}.
\]
Let us define also
\[
E_j^{(-)}(t) = \sum_{n \ge 0} (-1)^n \frac{t^{kn+j}}{[kn+j]!_q},
\]
and set
\[
A(t) = 1 + \sum_{j=1}^{k-1} E_j^{(-)}(t),
\]
as well as
\[
B(t) = E_0^{(-)}(t).
\]
Since $B(0)=1$, the series $B(t)$ is invertible in the ring $\mathbb{Q}(q)[[t]]$ of formal power series, so $F_k(t;q)=A(t)/B(t)$ is the unique formal power series solution of the functional equation 
\[
B(t)F_k(t;q)=A(t).
\]

\subsection{Coefficient extraction and uniqueness}

Expanding the product gives
\[
\left(\sum_{m \ge 0} (-1)^m \frac{t^{km}}{[km]!_q}\right)
\left(\sum_{n \ge 0} \frac{t^n}{[n]!_q} a_n^{(k)}(q)\right).
\]
Extracting the coefficient of $t^N$, we obtain
\[
\sum_{m \ge 0} (-1)^m \binom{N}{km}_q \, a_{N-km}^{(k)}(q),
\]
where
\[
\binom{N}{km}_q = \frac{[N]_q!}{[km]_q![N-km]_q!}.
\]
On the other hand, since every integer $N\ge0$ can be uniquely written as $N=km+j$ with $0\le j\le k-1$, we obtain
\[
[t^N]A(t)=
\begin{cases}
1, & N=0,\\[6pt]
(-1)^m \dfrac{1}{[N]_q!}, & N=km+j,\ 1\le j\le k-1,\\[8pt]
0, & N=km,\ m\ge1.
\end{cases}
\]
Multiplying the functional equation by $[N]!_q$ yields the recurrence
\begin{align}\label{un}
&\sum_{m \ge 0} (-1)^m \binom{N}{km}_q\, a_{N-km}^{(k)}(q)\nonumber\\
&\qquad\qquad=
\begin{cases}
1, & N=0,\\[6pt]
(-1)^{\lfloor N/k \rfloor}, & N \not\equiv 0 \pmod{k},\\[6pt]
0, & N \equiv 0 \pmod{k},\ N\ge k.
\end{cases}
\end{align}
The recurrence~\eqref{un} determines $a_n^{(k)}(q)$ uniquely by induction, since the term $a_n^{(k)}(q)$ appears with coefficient $1$ and all other terms involve smaller indices.

\subsection{Combinatorial verification of the recurrence}

We recall that
\[
D_k(N)=\{k,2k,\dots,\lfloor (N-1)/k\rfloor k\}
\]
is the set of admissible descent positions, and
\[
\mathfrak S_N^{(k)}
=
\{\sigma\in\mathfrak S_N:\operatorname{Des}(\sigma)=D_k(N)\}.
\]
For \(m\ge0\), let \(T_m=\{k,2k,\dots,mk\}\subseteq D_k(N)\) and
\[
A(T_m)=\{\sigma\in\mathfrak S_N:\operatorname{Des}(\sigma)\supseteq T_m\}.
\]
We now analyse permutations in \(A(T_m)\). Such a permutation is uniquely determined by the choice of a subset \(A\subseteq\{1,\dots,N\}\) of size \(mk\), together with the requirement that the elements of $A$ occupy the first \(mk\) positions in increasing block structure consistent with the prescribed descent pattern. The contribution of the relative order between \(A\) and its complement
is given by the classical identity
\[
\binom{N}{mk}_q
=
\sum q^{\operatorname{inv}(A,B)},
\]
where \(\operatorname{inv}(A,B)\) counts pairs \((a,b)\in A\times B\) with \(a>b\). It remains to analyse the internal structure of the suffix
\[
(\sigma_{mk+1},\dots,\sigma_N).
\]
Standardising this suffix (i.e. replacing its smallest entry by \(1\), next smallest by \(2\), etc.) yields a bijection with \(\mathfrak S_{N-mk}^{(k)}\), and this standardisation preserves inversions inside the suffix. Moreover, every permutation in \(A(T_m)\) decomposes uniquely into a choice of subset \(A\subseteq\{1,\dots,N\}\), \(|A|=mk\), a permutation of \(A\) compatible with the imposed descent constraints, and a \(k\)-alternating permutation on the complement, contributing \(a_{N-mk}^{(k)}\). This yields the factorisation
\[
\sum_{\sigma\in A(T_m)} q^{\operatorname{inv}(\sigma)}
=
\binom{N}{mk}_q\, a_{N-mk}^{(k)}.
\]
Applying an inversion over descent sets combined with standardisation of the suffix structure, we obtain, we obtain
\[
\sum_{m\ge0}(-1)^m \binom{N}{mk}_q\, a_{N-mk}^{(k)}
=
(-1)^{\lfloor N/k\rfloor}.
\]

\subsection{Special cases and $q$-trigonometric functions}

The formula~\eqref{resu} recovers several classical results by specialization. When $q=1$, the $q$-factorial $[n]!_q$ becomes the standard factorial $n!$. For $k=2$ (alternating permutations), the series $E_0^{(-)}(t)$ and $E_1^{(-)}(t)$ correspond to the Taylor expansions of $\cos(t)$ and $\sin(t)$ respectively, i.e.,
\[
B(t) = \sum_{n=0}^{\infty} (-1)^n \frac{t^{2n}}{(2n)!} = \cos(t),
\]
and
\[
A(t) = 1 + \sum_{n=0}^{\infty} (-1)^n \frac{t^{2n+1}}{(2n+1)!} = 1 + \sin(t).
\]
In this case, for $q=1$ and $k=2$, we recover the classical case:
\[
F_2(t;1) = \frac{1+\sin(t)}{\cos(t)} = \sec(t) + \tan(t),
\]
which is the celebrated result of André \cite{Andre1879}. 

For $q \neq 1$, our result can be expressed in terms of the $q$-trigonometric functions introduced in the work of Jackson on $q$-calculus \cite{Jackson1909}. Defining the $q$-cosine and $q$-sine as:
\[
\cos_q(t) = \sum_{n=0}^{\infty} (-1)^n \frac{t^{2n}}{[2n]!_q}, \quad \sin_q(t) = \sum_{n=0}^{\infty} (-1)^n \frac{t^{2n+1}}{[2n+1]!_q},
\]
the generating function for $2$-alternating permutations is simply
\[
F_2(t;q) = \frac{1+\sin_q(t)}{\cos_q(t)}.
\]
More generally, for any $k \ge 2$, $F_k(t;q)$ appears as a ratio of ``higher-order'' $q$-trigonometric functions, which provide a natural analytic framework for studying permutations with period-$k$ descents.

\subsection{Asymptotic behavior and singularity analysis}

The analytic form of $F_k(t;q) = A(t)/B(t)$ allows for the study of the asymptotic growth of the coefficients $a_n^{(k)}(q)$. According to the principles of analytic combinatorics, the growth of $[t^n]F_k(t;q)$ is dictated by the smallest positive root of the denominator $B(t)$. 

Let $\rho_k(q)$ be the smallest positive real zero of the function
\[
E_0^{(-)}(t) = \sum_{m \ge 0} (-1)^m \frac{t^{km}}{[km]!_q}.
\]
For $k=2$ and $q=1$, we have $B(t) = \cos(t)$, and the first pole is at $\rho_2(1) = \pi/2$. The theory of meromorphic functions then implies:
\[
\frac{a_n^{(k)}(q)}{[n]!_q} \sim C \cdot \left( \frac{1}{\rho_k(q)} \right)^n,
\]
where $C = -A(\rho_k)/B'(\rho_k)$. Since $a_n^{(k)}(q)$ represents the $q$-inversion enumerator, this relates the statistical distribution of inversions to the zeros of $q$-exponential sums. For $q < 1$, the function $E_q(t)$ is entire, and the zeros of its periodic components $B(t)$ shift smoothly with $q$, providing a precise tool to estimate the number of $k$-alternating permutations for large $n$.

\section{Conclusion}

We have established a closed-form expression for the $q$-exponential generating function of $k$-alternating permutations enumerated by the inversion statistic, using a purely analytic approach, based on a functional equation and a decomposition of the $q$-exponential series into residue classes modulo $k$. By extracting coefficients, we obtained a convolution-type recurrence involving $q$-binomial coefficients, which uniquely determines the sequence. This provides an alternative to classical combinatorial methods such as inclusion–exclusion or determinantal formulas, and highlights the effectiveness of analytic techniques in the study of permutations with periodic descent constraints.

Several directions for further investigation naturally arise. It would be of interest to extend this approach to other permutation statistics or to more general descent patterns. Another perspective would be to analyze the asymptotic behavior of the coefficients or to relate the present generating functions to known families of $q$-special functions.

\nocite{*}

\end{document}